\documentclass[11pt]{article}
\usepackage{mycommands}

\newcommand{\prefix}[1]{#1\mbox{-}\nobreak\hspace{0pt}}

\title{\textbf{The permutation classes \\
 $\av(\mathbf{1234},\mathbf{2341})$ and $\av(\mathbf{1243},\mathbf{2314})$}}

\author{$\phantom{{}^\dagger}$David Bevan${}^\dagger$}

\date{}

\begin{document}
\maketitle

{\let\thefootnote\relax\footnotetext
{${}^\dagger$Department of Mathematics and Statistics, The Open University, Milton Keynes, England.}}

{\let\thefootnote\relax}

\begin{abstract}
\noindent
We investigate the structure of the two permutation classes
defined by the sets of forbidden patterns $\{\mathbf{1234},\mathbf{2341}\}$ and $\{\mathbf{1243},\mathbf{2314}\}$.
By considering how the Hasse graphs of permutations in these classes can be built from a sequence of rooted source graphs,
we determine their algebraic generating functions.
Our approach is similar to that of ``adding a slice'', used previously to enumerate various classes of polyominoes and other combinatorial structures. To solve the relevant functional equations, we make extensive use of the kernel method.
\end{abstract}

\section{Introduction}
We consider a permutation to be simply an arrangement
of the numbers $1,2,\ldots n$ for some positive $n$.
A permutation $\pi$ is said to be \emph{contained} in another
permutation $\sigma$ if $\sigma$ has a subsequence whose terms
have the same relative ordering
as
those of $\pi$.
For example, $\mathbf{3241}$ is contained in $\mathbf{1573462}$ because the subsequence $\mathbf{5362}$ is ordered
in the same way as $\mathbf{3241}$. If $\pi$ is not contained in $\sigma$ then we say that $\sigma$ \emph{avoids} $\pi$.
For example, $\mathbf{1573462}$ avoids $\mathbf{3214}$.
In the context of containment and avoidance, a permutation is often called a \emph{pattern}.

The containment relation is a partial order on the set of all permutations, and a set of permutations closed downwards
(a down-set) in this partial order is called a \emph{permutation class}.
It is natural to define a permutation class by the minimal set of permutations that it avoids.
This minimal forbidden set of patterns is known as the \emph{basis} of the class.
The class with basis $B$ is denoted $\av(B)$.

Given a permutation class $\CCC$, we
denote by $\CCC_n$ the set of permutations in $\CCC$ of length $n$. The (univariate) \emph{generating function} of $\CCC$ is then
$\sum_{n\geqslant1}|\CCC_n|z^n = \sum_{\sigma\in\CCC}z^{|\sigma|}$, where $|\sigma|$ is the length of $\sigma$.
The \emph{growth rate} of $\CCC$ is defined by the limit
$
\liminfty\sqrt[n]{|\CCC_n|},
$
if it exists. It is widely believed (see the first conjecture in~\cite{Vatter2014}) that all permutation classes have a growth rate.

In the study of permutation classes, there has been significant interest in deriving the generating functions for classes with a few small basis elements (see~\cite{WikiEnumPermClasses} for an up-to-date list of results).
This has led to the enrichment of the theory of permutation classes due to the requisite development of a variety of enumeration techniques.
We add to this work by proving the following two theorems:

\thmbox{
\begin{thmO}\label{thmF}
The class of permutations avoiding $\mathbf{1234}$ and $\mathbf{2341}$ has the algebraic generating function
$$
\frac
{2-10\+z+9\+z^2+7\+z^3-4\+z^4 \:-\: (2-8\+z+9\+z^2-3\+z^3)\+\sqrt{1-4\+z}}
{(1-3\+z+z^2)\+\big((1-5\+z+4\+z^2) \:+\: (1-3\+z)\+\sqrt{1-4\+z}\big)} .
$$
Its growth rate is equal to 4.
\end{thmO}
}

\thmbox{
\begin{thmO}\label{thmE}
The class of permutations avoiding $\mathbf{1243}$ and $\mathbf{2314}$ has an algebraic generating function $F(z)$ which satisfies the cubic polynomial equation
$$
      (z-3\+z^2+2\+z^3)
\:-\: (1-5\+z+8\+z^2-5\+z^3)\+F(z)
\:+\: (2\+z-5\+z^2+4\+z^3)\+F(z)^2
\:+\: z^3\+F(z)^3
\;=\; 0
.
$$
Its growth rate is approximately 5.1955, the greatest real root of the quintic polynomial
$$
2 - 41\+z + 101\+z^2 - 97\+z^3 + 36\+z^4 -4\+z^5 .
$$
\end{thmO}
}

\vspace{9pt}
\textbf{Hasse graphs}

Corresponding to each permutation $\sigma$, we define an ordered plane graph $H_\sigma$, which we call its \emph{Hasse graph}.
If $P_\sigma$ is the poset on the points $(i,\sigma_i)$ in which $(i,\sigma_i)<(j,\sigma_j)$ if both $i<j$ and $\sigma_i<\sigma_j$, then $H_\sigma$ is the graph
corresponding to the Hasse diagram of $P_\sigma$.
See the figures throughout this paper for illustrations showing the Hasse graphs of permutations.
In practice, we tend not to distinguish between a permutation and its Hasse graph.
The minimal elements of the poset $P_\sigma$ are known as the \emph{left-to-right minima} of the permutation $\sigma$.
Similarly, maximal elements of $P_\sigma$ are called \emph{right-to-left maxima} of $\sigma$.

Hasse graphs of permutations were previously considered by Bousquet-M\'elou \& Butler~\cite{BMB2007}, who
determined the algebraic generating function of the family of \emph{forest-like} permutations
whose Hasse graphs are acyclic.
More recently, they
have been used by the present author~\cite{Bevan2014} to
establish
a new lower bound for the growth rate of $\av(\mathbf{1324})$.

\begin{figure}[t] 
  $$
  \begin{tikzpicture}[scale=0.225,line join=round]
    \draw [red,thin] (3,11)--(5,16);
    \draw [red,thin] (5,16)--(4,4)--(7,14);
    \draw [red,thin] (12,13)--(8,8)--(13,12)--(9,6)--(12,13);
    \draw [red,thin] (6,1)--(7,14);
    \draw [red,thin] (8,8)--(6,1)--(9,6);
    \draw [red,thin] (10,3)--(12,13)--(11,2)--(13,12)--(10,3)--(14,10)--(11,2)--(15,9)--(10,3)--(16,7)--(11,2)--(17,5)--(10,3);
    \draw [blue,very thick] (2,17)--(1,15)--(5,16);
    \draw [blue,very thick] (7,14)--(3,11);
    \draw [blue,very thick] (12,13)--(3,11)--(13,12);
    \draw [blue,very thick] (16,7)--(9,6)--(4,4)--(8,8)--(14,10)--(9,6)--(15,9)--(8,8);
    \draw [blue,very thick] (4,4)--(17,5);
    \draw [blue,very thick] (10,3)--(6,1)--(11,2);
    \plotpermnobox{}{15,17,11, 4,16, 1,14, 8, 6, 3, 2,13,12,10, 9, 7, 5}
    \draw [thin] (1,15)  circle [radius=0.4];
    \draw [thin] (3,11) circle [radius=0.4];
    \draw [thin] (4,4) circle [radius=0.4];
    \draw [thin] (6,1) circle [radius=0.4];
  \end{tikzpicture}
  $$
  \caption{
           A permutation in class $\FFF$, spanned by four source graphs} 
  \label{figFPerm}
\end{figure}
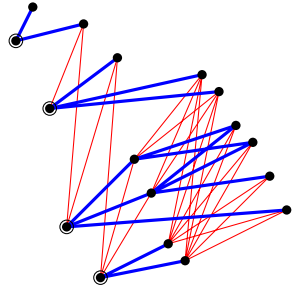

Given a permutation $\sigma$, we partition the vertices of $H_\sigma$ by spanning it with a sequence of graphs, which we call the \emph{source graphs} of $\sigma$.
There is one source graph for each left-to-right minimum of $\sigma$.
Suppose $u_1,\ldots,u_m$ are the vertices of $H_\sigma$ corresponding to the left-to-right minima of $\sigma$, listed from left to right.
Then the $k$th source graph $G_k$ is the graph induced by $u_k$ and those vertices of $H_\sigma$ lying above and to the right of $u_k$ that are not in $G_1,\ldots,G_{k-1}$. We refer to $u_k$ as the \emph{root} of source graph $G_k$.
See Figure~\ref{figFPerm} for an illustration.
The structure of the source graphs of permutations in a specific permutation class is constrained by the need to avoid the patterns in the basis of the class.
If the source graphs for some class are acyclic, we refer to them as \emph{source trees}.

The \emph{bottom subgraph}
of a Hasse graph
is
the graph induced by its lowest vertex (the
least entry in the permutation) and all the vertices lying above and to its right.
Observe that the bottom subgraph may contain vertices from more than one source graph.
For example, the bottom subgraph of the Hasse graph in Figure~\ref{figFPerm} contains vertices from three source graphs.
Bottom subgraphs of permutations in a specific permutation class satisfy the same structural restrictions as do the source graphs.
We refer to an acyclic bottom subgraph as a \emph{bottom subtree}.

We build the Hasse graph of a permutation by starting with a source graph and then repeatedly adding another source graph to the lower right.
The technique is similar to that of ``adding a slice'', which has been used to enumerate constrained compositions and other classes of polyominoes, a topic of interest in statistical mechanics
(see, for example, Bousquet-M\'elou's review paper~\cite{Bousquet-Melou1996}, the books of van Rensburg~\cite{vanRensburg2000} and Guttmann~\cite{Guttmann2009}, and Flajolet \& Sedge\-wick~\cite[Examples~III.22 and~V.20]{FS2009}).
When a source graph is added, its vertices are interleaved horizontally with the non-root vertices of the bottom subgraph of the graph built from the previous source graphs.
Typically, the positioning of the vertices of the new source graph is constrained by the need to avoid forbidden patterns.

In order to derive the univariate generating functions we require,
we make use of multivariate functions involving additional ``catalytic'' variables that record certain parameters of the bottom subgraph of the permutations.
These additional variables enable us to establish recurrence relations which we can then solve using the kernel method.
Typically, when employing a multivariate generating function, we treat it simply as a function of the relevant catalytic variable, writing, for example, $F(u)$ rather than $F(z,u)$.

Occasionally, we also
make use of a variant of the symbolic structural notation presented in Flajolet \& Sedge\-wick~\cite{FS2009} 
to establish functional equations.
In particular,
$\ZZZ$ is the atomic class consisting of a single vertex, and
we use
$\seq{\AAA}$ to represent a possibly empty sequence of elements of $\AAA$ and
$\seqplus{\AAA}$ to represent a non-empty sequence of elements of $\AAA$.

The two classes we enumerate are quite distinct structurally.
A source graph in
$\av(\mathbf{1234},\mathbf{2341})$
consists of a root together with a \prefix{$\mathbf{123}$}avoider formed from the non-root vertices.
However, the presence of a $\mathbf{123}$ forces any subsequent source graph to be simply a fan.
In contrast, $\av(\mathbf{1243},\mathbf{2314})$ has plane source graphs and a much more uniform structure.
We enumerate
$\av(\mathbf{1234},\mathbf{2341})$ in Section~\ref{sectF}.
In doing so, the kernel method is used six times to solve
the relevant functional equations.
The class $\av(\mathbf{1243},\mathbf{2314})$ is enumerated in
Section~\ref{sectE}. This requires an
unusual simultaneous double application of the kernel method.

\section{Permutations avoiding \texorpdfstring{$\mathbf{1234}$ and $\mathbf{2341}$}{1234 and 2341}}\label{sectF}

Let us use $\FFF$ to denote $\av(\mathbf{1234},\mathbf{2341})$.
The structure of class $\FFF$ depends critically on the presence or absence of occurrences of the pattern $\mathbf{123}$.
In light of this,
to enumerate this class, we partition it into three sets $\AAA$, $\BBB$ and $\CCC$ as follows:
\begin{bullets}
  \item $\AAA=\av(\mathbf{123})$.
  \item $\BBB=\av(\mathbf{1234},\mathbf{2341},\mathbf{13524},\mathbf{14523})\setminus\AAA$. Every permutation in $\BBB$ contains at least one occurrence of a $\mathbf{123}$, but avoids $\mathbf{13524}$ and $\mathbf{14523}$.
  \item $\CCC=\av(\mathbf{1234},\mathbf{2341})\setminus(\AAA\cup\BBB)$. Every permutation in $\CCC$ contains a $\mathbf{13524}$ or a $\mathbf{14523}$.
\end{bullets}
We refer to a permutation in $\AAA$ as an \prefix{$\AAA$}permutation, and similarly for $\BBB$ and $\CCC$.

The addition of a source graph to a \prefix{$\CCC$}permutation can only yield another \prefix{$\CCC$}permutation (since it can't cause the removal a $\mathbf{13524}$ or $\mathbf{14523}$ pattern).
Similarly, the addition of a source graph to a \prefix{$\BBB$}permutation can't result in an \prefix{$\AAA$}permutation.
Hence, we can enumerate $\AAA$ without first considering $\BBB$ or $\CCC$, and can enumerate $\BBB$ before considering $\CCC$.

Before investigating the structure of permutations in $\AAA$, $\BBB$ and $\CCC$, let us briefly examine what a typical source graph in $\FFF$ looks like. Firstly, the avoidance of $\mathbf{1234}$ means that the non-root vertices of any source graph form a \prefix{$\mathbf{123}$}avoider. 
Secondly,
the avoidance of $\mathbf{2341}$ presents no additional restriction on the structure of a source graph,
because the presence of a $\mathbf{2341}$ would
imply the presence of a $\mathbf{123}$ in the non-root vertices.
Thus a source graph in $\FFF$ consists of a root together with a \prefix{$\mathbf{123}$}avoider formed from the non-root vertices.

\subsubsection*{The structure of set $\AAA$}

We begin by looking at $\AAA=\av(\mathbf{123})$.
As is very well known, this class is enumerated by the Catalan numbers.
However, we need to keep track of the structure of the bottom graph. So we must determine the appropriate
bivariate generating function.

Let $\AAA_\ssS$ denote the set of source graphs in set $\AAA$.
Now, each member of $\AAA_\ssS$ is simply a \emph{fan}, a root vertex connected to a (possibly empty) sequence of pendant edges.
Bottom subgraphs are also fans.
Thus source graphs and bottom subgraphs of $\AAA$ are acyclic.

When enumerating $\AAA$, we use the variable $u$ to mark the \emph{number of leaves} (non-root vertices) in the bottom subtree.
The generating function for $\AAA_\ssS$ is thus given by
$$
A_\ssS(u) \;=\; z+z^2\+u+z^3\+u^2+\ldots \;=\; \frac{z}{1-z\+u} .
$$

We now consider the process of building an \prefix{$\AAA$}permutation from a sequence of source trees.
When a source tree is added to an \prefix{$\AAA$}permutation,
the root vertex of the source tree may be inserted to the left of zero or more of the leaves of the bottom subtree.
See Figure~\ref{figABuild} for an illustration.
Note that, in this and other similar figures, the original bottom subgraph is displayed to the upper left, with the new source graph to the lower right.

The action of adding a source tree is thus seen to be reflected by the linear operator $\oper_{\ssA}$ whose effect on $u^k$ is given by
$$
\oper_{\ssA}\big[u^k\big] \;=\; A_\ssS(u)\+(1+u+\ldots+u^k) \;=\; A_\ssS(u)\+\frac{1-u^{k+1}}{1-u}.
$$
Hence,
the bivariate generating function $A(u)$ for $\AAA$ is defined by the following recursive functional equation:
\begin{equation*}
   A(u)  \;=\;  A_\ssS(u) \:+\: A_\ssS(u)\+\frac{A(1)-u\+A(u)}{1-u}.
\end{equation*}
\begin{figure}[t] 
  $$
  \begin{tikzpicture}[scale=0.225,line join=round]
    \draw [very thick] (2,11)--(1,6)--(3,10);
    \draw [very thick] (5,9)--(1,6)--(6,8);
    \draw [very thick] (7,7)--(1,6);
    \draw [blue,very thick] (8,5)--(4,1)--(9,4);
    \draw [blue,very thick] (10,3)--(4,1)--(11,2);
    \draw [red,thin] (5,9)--(4,1)--(6,8);
    \draw [red,thin] (7,7)--(4,1);
    \plotpermnobox{}{6,11,10, 1, 9, 8, 7, 5, 4, 3, 2}
    \draw [thin] (1,6)  circle [radius=0.4];
    \draw [thin] (4,1) circle [radius=0.4];
  \end{tikzpicture}
  $$
  \caption{Adding a source tree to the bottom subtree of an \prefix{$\AAA$}permutation}
  \label{figABuild}
\end{figure}
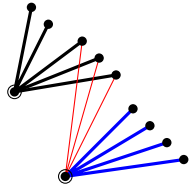This equation can be solved using the \emph{kernel method}.
To start, we express $A(u)$ in terms of $A(1)$,
by expanding and rearranging to give
\begin{equation}\label{eqAKernel}
A(u)  \;=\;   \frac{z \+\big(1-u+A(1)\big)}{1-u+z\+u^2}.
\end{equation}
Equivalently, we have the equation
$$
(1-u+z\+u^2)\+A(u)  \;=\;   z \+\big(1-u+A(1)\big) .
$$
Now, if we set $u$ to be a root of the multiplier of $A(u)$, we obtain a linear equation for $A(1)$.
This is known as ``cancelling the kernel'' (the multiplier of $A(u)$ being the kernel).
The appropriate root to use can be identified from
the combinatorial requirement that the series expansion of $A(1)$ contains no negative exponents and
has only non-negative coefficients.

In this case, the correct root is $u=(1-\sqrt{1-4\+z})/2\+z$, which yields
the univariate generating function for $\AAA$,
$$
A(1) \;=\; \frac{1-\sqrt{1-4\+ z}}{2\+ z}-1.
$$
This is the generating function for the Catalan numbers as expected.

Finally, by substituting for $A(1)$ back into~\eqref{eqAKernel} we
get the following explicit algebraic expression for $A(u)$:
$$
A(u)  \;=\; \frac{1-2\+z\+u-\sqrt{1-4\+z}}{2\+(1-u+z\+u^2)} .
$$
On this occasion, we have explained every step of the derivation. On subsequent occasions, we present fewer details of the algebraic manipulations. 

\subsubsection*{The structure of set $\BBB$}

We now consider set $\BBB$.
Recall that sets $\BBB$ and $\CCC$ consist of those permutations in class $\FFF$ that contain at least one occurrence of a $\mathbf{123}$.
We need to keep track of the position of the
leftmost occurrence of a $\mathbf{3}$ in such a pattern.
Given a permutation in $\BBB$ or $\CCC$, let us call the vertex corresponding to the leftmost $\mathbf{3}$ in a $\mathbf{123}$ the \emph{spike}.
In the figures, the spike is marked with a star.

We now make a key observation.\label{obsKey}
When adding a source graph to a permutation containing a $\mathbf{123}$, no vertex of the source graph may be positioned to the right of the spike, or else a $\mathbf{2341}$ would be created.
Hence, the spike in any permutation in classes $\BBB$ or $\CCC$ occurs in its bottom subgraph.
When enumerating sets $\BBB$ and $\CCC$, we use the variable $u$ to mark \emph{the
number of vertices to the left of the spike} in its bottom subgraph.

\begin{figure}[ht]
  $$
  \begin{tikzpicture}[scale=0.225,line join=round]
    \draw [blue,very thick] (2,15)--(1,1)--(3,14);
    \draw [blue,very thick] (4,13)--(1,1)--(5,10);
    \draw [blue,very thick] (6,9)--(1,1)--(7,7);
    \draw [blue,very thick] (8,4)--(1,1);
    \draw [blue,very thick] (14,3)--(1,1)--(15,2);
    \draw [blue,very thick] (9,12)--(5,10)--(10,11)--(6,9)--(9,12)--(7,7)--(10,11)--(8,4)--(9,12);
    \draw [blue,very thick] (7,7)--(11,8)--(8,4);
    \draw [blue,very thick] (12,6)--(8,4)--(13,5);
    \plotpermnobox{}{1,15,14,13,10,9,7,4,12,11, 8, 6, 5, 3, 2}
    \draw [thin] (1,1)  circle [radius=0.4];
    \node[above] at (9,11.6) {${}^\bigstar$};
  \end{tikzpicture}
  $$
  \caption{A source graph in set $\BBB$} 
  \label{figBSource}
\end{figure}
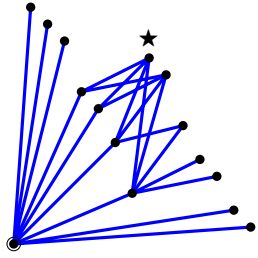
Let $\BBB_\ssS$ be the set of source graphs in set $\BBB$.
Since \prefix{$\BBB$}permutations contain a $\mathbf{123}$ but
avoid $\mathbf{13524}$ and $\mathbf{14523}$,
the non-root vertices of a permutation in $\BBB_\ssS$ consist of two descending sequences, the second sequence beginning (with the spike) above the last vertex in the first sequence.
See Figure~\ref{figBSource} for an illustration.
If we consider the non-root vertices in order from top to bottom, then it can be seen that
$\BBB_\ssS$ is defined by the structural equation
$$
\BBB_\ssS \;=\; u\+\ZZZ \:\times\: \seq{u\+\ZZZ} \:\times\: \ZZZ \:\times\: \seq{u\+\ZZZ+\ZZZ} \:\times\: u\+\ZZZ \:\times\: \seq{\ZZZ}.
$$
The first term on the right corresponds to the root and the remaining terms deal with the non-root vertices in order from top to bottom, vertices to the left of the spike being marked with $u$.
The third term corresponds to the spike and the fifth represents the lowest point to the left of the spike (the rightmost $\mathbf{2}$ of a $\mathbf{123}$).
Hence, the generating function for $\BBB_\ssS$ is
$$
B_\ssS(u) \;=\; \frac{z^3\+u^2}{(1-z)\+ (1-z\+u)\+ (1-z-z\+u)} .
$$

We now study the process of building a \prefix{$\BBB$}permutation from a sequence of source graphs.
There are two cases. A permutation in $\BBB$ may result either from the addition of a source graph to an \prefix{$\AAA$}permutation, or else
from adding a source graph to another \prefix{$\BBB$}permutation. We address these two cases in turn.

\begin{figure}[ht]
  $$
  \begin{tikzpicture}[scale=0.225,line join=round]
    \draw [very thick] (2,17)--(1,12)--(4,16);
    \draw [very thick] (5,15)--(1,12)--(6,14);
    \draw [very thick] (7,13)--(1,12);
    \draw [blue,very thick] (8,11)--(3,2)--(9,9);
    \draw [blue,very thick] (10,8)--(3,2)--(11,5);
    \draw [blue,very thick] (12,4)--(3,2);
    \draw [blue,very thick] (16,3)--(3,2);
    \draw [blue,very thick] (9,9)--(13,10)--(10,8);
    \draw [blue,very thick] (11,5)--(13,10)--(12,4)--(14,7)--(11,5)--(15,6)--(12,4);
    \draw [red,thin] (5,15)--(3,2)--(4,16);
    \draw [red,thin] (7,13)--(3,2)--(6,14);
    \plotpermnobox{}{12,17, 2,16,15,14,13,11, 9, 8, 5, 4,10, 7, 6, 3}
    \draw [thin] (1,12)  circle [radius=0.4];
    \draw [thin] (3,2) circle [radius=0.4];
    \node[above] at (13,9.6) {${}^\bigstar$};
  \end{tikzpicture}
  \qquad\qquad
  \begin{tikzpicture}[scale=0.225,line join=round]
    \draw [very thick] (2,17)--(1,12)--(4,16);
    \draw [very thick] (5,15)--(1,12)--(11,14);
    \draw [very thick] (12,13)--(1,12);
    \draw [blue,very thick] (6,11)--(3,2)--(7,10);
    \draw [blue,very thick] (8,9)--(3,2)--(9,8);
    \draw [blue,very thick] (10,7)--(3,2)--(13,6);
    \draw [blue,very thick] (14,5)--(3,2)--(15,4);
    \draw [blue,very thick] (16,3)--(3,2);
    \draw [red,thin] (5,15)--(3,2)--(4,16);
    \draw [red,thin] (11,14)--(6,11)--(12,13)--(7,10)--(11,14)--(8,9)--(12,13)--(9,8)--(11,14)--(10,7)--(12,13);
    \plotpermnobox{}{12,17, 2,16,15,11,10, 9, 8, 7,14,13, 6, 5, 4, 3}
    \draw [thin] (1,12)  circle [radius=0.4];
    \draw [thin] (3,2) circle [radius=0.4];
    \node[above] at (11,13.6) {${}^\bigstar$};
  \end{tikzpicture}
  \qquad\qquad
  \begin{tikzpicture}[scale=0.225,line join=round]
    \draw [very thick] (2,17)--(1,12)--(4,16);
    \draw [very thick] (5,15)--(1,12)--(11,14);
    \draw [very thick] (12,13)--(1,12);
    \draw [blue,very thick] (6,11)--(3,2)--(7,9);
    \draw [blue,very thick] (8,8)--(3,2)--(9,5);
    \draw [blue,very thick] (10,4)--(3,2);
    \draw [blue,very thick] (16,3)--(3,2);
    \draw [blue,very thick] (7,9)--(13,10)--(8,8);
    \draw [blue,very thick] (9,5)--(13,10)--(10,4)--(14,7)--(9,5)--(15,6)--(10,4);
    \draw [red,thin] (5,15)--(3,2)--(4,16);
    \draw [red,thin] (11,14)--(6,11)--(12,13)--(7,9)--(11,14)--(8,8)--(12,13)--(9,5)--(11,14)--(10,4)--(12,13);
    \plotpermnobox{}{12,17, 2,16,15,11, 9, 8, 5, 4,14,13,10, 7, 6, 3}
    \draw [thin] (1,12)  circle [radius=0.4];
    \draw [thin] (3,2) circle [radius=0.4];
    \node[above] at (11,13.6) {${}^\bigstar$};
  \end{tikzpicture}
  $$
  \caption{Ways to create a \prefix{$\BBB$}permutation by adding a source graph to the bottom subtree of an \prefix{$\AAA$}permutation}
  \label{figABBuild}
\end{figure}
One way to create a \prefix{$\BBB$}permutation from an \prefix{$\AAA$}permutation is to
add a source graph from $\BBB_\ssS$, positioning its root to the left of zero or more of the leaves of the bottom subtree
of the \prefix{$\AAA$}permutation
and its non-root vertices to the right of the bottom subtree.
In this case, the new permutation inherits its spike from the added source graph.
This is illustrated in the left diagram in Figure~\ref{figABBuild}.
The generating function for this set of permutations is thus given by
$$
B_{\textsf{AB1}}(u) \;=\; B_\ssS(u)\+\frac{A(1)-u\+A(u)}{1-u}
.
$$
For simplicity, we choose not to present the expanded form of $B_{\textsf{AB1}}(u)$, or that of most subsequent expressions. They can all be represented in the form $(p+q\+\sqrt{1-4\+z})/r$ for appropriate polynomials $p$, $q$ and $r$.

The other possibility for creating a \prefix{$\BBB$}permutation from an \prefix{$\AAA$}permutation involves the positioning of some non-root vertices of the source graph to the left of some of the leaves in the bottom subtree, making one of the vertices in the original bottom subtree the spike.
The source graph may be drawn from either $\AAA_\ssS$ or $\BBB_\ssS$, as illustrated in the centre and right diagrams in Figure~\ref{figABBuild}.

In this situation, if the source graph has a spike,
it must be positioned to the right of all leaves in the bottom subtree, or else a $\mathbf{1234}$ would be created.
Furthermore, any source graph vertices placed to the left of leaves in the bottom subtree
must occur at the same position in the bottom subtree, or else a $\mathbf{13524}$ would be created.
This position may be chosen independently of where the root vertex is placed.

From these considerations, it can be determined that
the resulting set of permutations has the generating function defined by
$$
B_{\textsf{AB2}}(u) \;=\; \Big(B_\ssS(u)+\frac{z^2\+u^2}{(1-z)\+(1-z\+u)}\Big) \+ \frac{1}{1-u} \+ \Big(A\!'(1)-\frac{u}{1-u}\+\big(A(1)-A(u)\big)\Big) ,
$$
where the presence of the derivative $A\!'$ is a consequence of the independent choice of two positions in the bottom tree.


\begin{figure}[ht]
  $$
  \begin{tikzpicture}[scale=0.225,line join=round]
    \draw [very thick] (2,12)--(1,5)--(3,10)--(10,11);
    \draw [very thick] (5,9)--(1,5)--(6,7);
    \draw [very thick] (1,5)--(12,6);
    \draw [very thick] (11,8)--(6,7)--(10,11)--(5,9);
    \draw [blue,very thick] (7,4)--(4,1)--(8,3);
    \draw [blue,very thick] (9,2)--(4,1);
    \draw [red,thin] (5,9)--(4,1)--(6,7);
    \draw [red,thin] (10,11)--(7,4)--(11,8)--(8,3)--(10,11)--(9,2)--(11,8);
    \draw [red,thin] (7,4)--(12,6)--(8,3);
    \draw [red,thin] (9,2)--(12,6);
    \plotpermnobox{}{5,12,10, 1, 9, 7, 4, 3, 2,11, 8, 6}
    \draw [thin] (1,5)  circle [radius=0.4];
    \draw [thin] (4,1) circle [radius=0.4];
    \node[above] at (10,10.6) {${}^\bigstar$};
  \end{tikzpicture}
  $$
  \caption{Adding a source tree to the bottom subgraph of a \prefix{$\BBB$}permutation}
  \label{figBBBuild}
\end{figure}
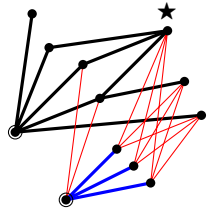
Finally, we consider the addition of a source graph to a \prefix{$\BBB$}permutation.
As we noted in our key observation on page~\pageref{obsKey},
no vertex of the source graph may be positioned to the right of the spike in the bottom subgraph.
As a result, the new source graph may not contain a $\mathbf{123}$ or else a $\mathbf{1234}$ would be created, so the source graph must be a member of $\AAA_\ssS$ (a fan).
Moreover, the leaves of the source tree must be positioned \emph{immediately} to the left of the spike, or else a $\mathbf{1234}$ would be created.
See Figure~\ref{figBBBuild} for an illustration.

Note that, as a consequence of these restrictions, it is impossible for the addition of a source graph to a \prefix{$\BBB$}permutation to create a $\mathbf{13524}$ or $\mathbf{14523}$. So it is not possible to extend a \prefix{$\BBB$}permutation so as to create a \prefix{$\CCC$}permutation.

Thus
the bivariate generating function $B(u)$ of set $\BBB$ is defined by the following recursive functional equation:
\begin{equation*}
   B(u)  \;=\;  B_\ssS(u) + B_{\textsf{AB1}}(u) + B_{\textsf{AB2}}(u) \:+\: \frac{z\+u}{1-z\+u}\+\frac{B(1)-B(u)}{1-u} ,
\end{equation*}
where the final term reflects the addition of a source tree to a \prefix{$\BBB$}permutation.

This equation is amenable to the kernel method.
It can be rearranged to express $B(u)$ in terms of $B(1)$. The kernel can then be cancelled by setting $u=(1-\sqrt{1-4\+z})/2\+z$, which yields an expression for $B(1)$:
$$
B(1)
\;=\;
\frac{-1+8\+z-19\+z^2+12\+z^3 \:+\: (1-6\+z+9\+z^2-2\+z^3)\+\sqrt{1-4 z}}{2\+z^3\+(1-4\+z)}.
$$
Substitution then
results in an explicit algebraic expression for $B(u)$, which we refrain from presenting explicitly due to its size.

\vspace{9pt}
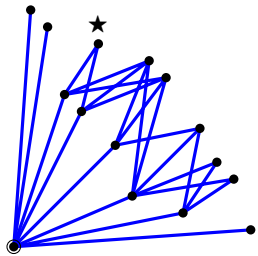
\begin{figure}[ht]
  $$
  \begin{tikzpicture}[scale=0.225,line join=round]
    \draw [blue,very thick] (2,15)--(1,1)--(3,14);
    \draw [blue,very thick] (1,1)--(4,10)--(6,13)--(5,9);
    \draw [blue,very thick] (5,9)--(1,1)--(7,7);
    \draw [blue,very thick] (8,4)--(1,1);
    \draw [blue,very thick] (12,8)--(11,3)--(1,1)--(15,2);
    \draw [blue,very thick] (9,12)--(4,10)--(10,11)--(5,9)--(9,12)--(7,7)--(10,11)--(8,4)--(9,12);
    \draw [blue,very thick] (7,7)--(12,8)--(8,4);
    \draw [blue,very thick] (13,6)--(8,4)--(14,5)--(11,3)--(13,6);
    \plotpermnobox{}{1,15,14,10, 9,13, 7, 4,12,11, 3, 8, 6, 5, 2}
    \draw [thin] (1,1)  circle [radius=0.4];
    \node[above] at (6,12.6) {${}^\bigstar$};
  \end{tikzpicture}
  $$
  \caption{A source graph in set $\CCC$}
  \label{figCSource}
\end{figure}
\subsubsection*{The structure of set $\CCC$}

We begin our enumeration of $\CCC$ by counting its set of source graphs, which we denote $\CCC_\ssS$.
Rather than doing this directly,
we enumerate all the source graphs that contain a $\mathbf{123}$ (i.e.~those in either $\BBB_\ssS$ or $\CCC_\ssS$) and then subtract those in $\BBB_\ssS$.
To begin, we consider how we might build an \emph{arbitrary} source graph in class $\FFF$ by adding vertices from left to right.

Suppose we have a partly formed source graph with at least one non-root vertex, whose rightmost vertex is $v$,
and we want to add further vertices to its right.
What are the options?
If $v$ is not the lowest of the non-root vertices, then any subsequent vertices must be placed lower than~$v$.
The only other restriction is that vertices must be positioned higher than the root.
If we use $y$ to mark the number of positions in which a vertex may be inserted, then the action of adding a new vertex can be seen to be reflected by the following linear operator: 
\begin{equation*}
  \oper_{\textsf{\L}}\big[f(y)\big] \;=\; z\+y^2\+\frac{f(1)-f(y)}{1-y}.
\end{equation*}
We choose to denote this operator $\oper_{\textsf{\L}}$ because it corresponds to the action used in building a {\L}uka\-sie\-wicz path.

Now let us consider source graphs that have no vertices to the right of the spike.
These are in~$\BBB_\ssS$, so let's call this set $\BBB_\textsf{S0}$.
As usual, we mark with $u$ those vertices to the left of the spike.
If, in addition, we mark with $y$ those vertices not above the spike, then $\BBB_\textsf{S0}$ is defined by the structural equation
$$
\BBB_\textsf{S0} \;=\; u\+\ZZZ \:\times\: \seq{u\+\ZZZ} \:\times\: \seqplus{u\+y\+\ZZZ} \:\times\: y\+\ZZZ .
$$
It is readily seen that $y$ correctly marks the number of positions in which an additional vertex may be inserted to the right.

Let $\DDD_\ssS=\BBB_\ssS\cup\CCC_\ssS$. Since every member of $\DDD_\ssS$ is built from an element of $\BBB_\textsf{S0}$ by applying $\oper_{\textsf{\L}}$ zero or more times, it follows that the generating function for $\DDD_\ssS$ is defined by the recursive functional equation
$$
D_\ssS(y) \;=\; \frac{z^3\+y^2\+u^2}{(1-z\+u)\+(1-z\+y\+u)} \:+\: z\+y^2\+\frac{D_\ssS(1)-D_\ssS(y)}{1-y} .
$$
This equation can be solved for $D_\ssS(1)$ by the kernel method,
using $y=(1-\sqrt{1-4\+z})/2\+z$ to
cancel the kernel.
The generating function for $\CCC_\ssS$ 
is then defined by
$$
C_\ssS(u) \;=\; D_\ssS(1) \:-\: B_\ssS(u) .
$$

We now study the process of building a \prefix{$\CCC$}permutation from a sequence of source graphs.
As with set $\BBB$, there are two cases. A permutation in $\CCC$ may result either from the addition of a source graph to an \prefix{$\AAA$}permutation, or else
from adding a source graph to another \prefix{$\CCC$}permutation.
(As we observed above, it is not possible to create a \prefix{$\CCC$}permutation by adding a source graph to a \prefix{$\BBB$}permutation.)
We address the two cases in turn.

\begin{figure}[ht]
  $$
  \begin{tikzpicture}[scale=0.225,line join=round]
    \draw [very thick] (2,19)--(1,14)--(4,18);
    \draw [very thick] (5,17)--(1,14)--(6,16);
    \draw [very thick] (7,15)--(1,14);
    \draw [blue,very thick] (8,13)--(3,1)--(9,10);
    \draw [blue,very thick] (10,9)--(3,1)--(11,8);
    \draw [blue,very thick] (12,7)--(3,1)--(13,5);
    \draw [blue,very thick] (15,4)--(3,1)--(16,3);
    \draw [blue,very thick] (3,1)--(19,2);
    \draw [blue,very thick] (14,12)--(10,9)--(17,11)--(9,10)--(14,12)--(13,5);
    \draw [blue,very thick] (11,8)--(14,12)--(12,7)--(17,11)--(11,8);
    \draw [blue,very thick] (17,11)--(13,5)--(18,6)--(15,4)--(17,11)--(16,3)--(18,6);
    \draw [red,thin] (4,18)--(3,1)--(5,17);
    \draw [red,thin] (6,16)--(3,1)--(7,15);
    \plotpermnobox{}{14,19, 1,18,17,16,15,13,10, 9, 8, 7, 5,12, 4, 3,11, 6, 2}
    \draw [thin] (1,14)  circle [radius=0.4];
    \draw [thin] (3,1) circle [radius=0.4];
    \node[above] at (14,11.6) {${}^\bigstar$};
  \end{tikzpicture}
  \qquad\qquad\qquad\quad
  \begin{tikzpicture}[scale=0.225,line join=round]
    \draw [very thick] (2,19)--(1,14)--(4,18);
    \draw [very thick] (8,17)--(1,14)--(10,16);
    \draw [very thick] (13,15)--(1,14);
    \draw [blue,very thick] (5,13)--(3,1)--(6,10);
    \draw [blue,very thick] (7,9)--(3,1)--(9,8);
    \draw [blue,very thick] (11,7)--(3,1)--(12,5);
    \draw [blue,very thick] (15,4)--(3,1)--(16,3);
    \draw [blue,very thick] (3,1)--(19,2);
    \draw [blue,very thick] (14,12)--(7,9)--(17,11)--(6,10)--(14,12)--(12,5);
    \draw [blue,very thick] (9,8)--(14,12)--(11,7)--(17,11)--(9,8);
    \draw [blue,very thick] (17,11)--(12,5)--(18,6)--(15,4)--(17,11)--(16,3)--(18,6);
    \draw [red,thin] (4,18)--(3,1);
    \draw [red,thin] (8,17)--(5,13)--(10,16)--(6,10)--(8,17)--(7,9)--(10,16)--(9,8)--(13,15);
    \draw [red,thin] (5,13)--(13,15)--(6,10);
    \draw [red,thin] (7,9)--(13,15)--(11,7);
    \draw [red,thin] (12,5)--(13,15);
    \plotpermnobox{}{14,19, 1,18,13,10, 9,17, 8,16, 7, 5,15,12, 4, 3,11, 6, 2}
    \draw [thin] (1,14)  circle [radius=0.4];
    \draw [thin] (3,1) circle [radius=0.4];
    \node[above] at (8,16.6) {${}^\bigstar$};
  \end{tikzpicture}
  $$
  \caption{Ways to create a \prefix{$\CCC$}permutation by adding a source graph to the bottom subtree of an \prefix{$\AAA$}permutation}
  \label{figACBuild}
\end{figure}
One way to create a \prefix{$\CCC$}permutation from an \prefix{$\AAA$}permutation is to
add a source graph from $\CCC_\ssS$, positioning its root to the left of zero or more of the leaves of the bottom subtree
of the \prefix{$\AAA$}permutation
and its non-root vertices to the right of the bottom subtree.
This is illustrated in the left diagram in Figure~\ref{figACBuild}.
The generating function for this set of permutations is thus
$$
C_{\textsf{AC1}}(u) \;=\; C_\ssS(u)\+\frac{A(1)-u\+A(u)}{1-u} .
$$

The other method for creating a \prefix{$\CCC$}permutation from an \prefix{$\AAA$}permutation involves the positioning of some non-root vertices of the source graph to the left of some of the leaves in the bottom subtree.
This is illustrated in the right diagram in Figure~\ref{figACBuild}.
In analysing this method, it is more convenient to look, more generally, at how an \prefix{$\AAA$}permutation can be extended to yield a permutation containing a $\mathbf{123}$, in either $\BBB$ or $\CCC$. We can then subtract those members of $\BBB$ that are enumerated by $B_{\textsf{AB2}}$.

We achieve the enumeration by adding vertices from left to right in four steps:
\begin{bulletnums}
\item The first step adds the root.
\item The second step adds the first non-root vertex, which determines the position of the new spike, and also any other vertices positioned to the left of the spike.
\item The third step adds any additional vertices to the right of the spike but to the left of some other leaves in the bottom subtree. The addition of such vertices creates occurrences of $\mathbf{13524}$.
\item Finally, the fourth step adds any vertices to the right of the bottom subtree.
\end{bulletnums}

Step 1: Permutations that result from the addition of the root vertex are enumerated by
$$
D_\textsf{1}(u) \;=\; z\+u\+\frac{A(1)-A(u)}{1-u} .
$$
Step 2: In this step, we insert the descending sequence of vertices that creates the new spike.
In the generating function for permutations resulting from this action, we introduce two additional catalytic variables that we require for steps 3 and 4.
For use in step 3, $v$ marks the number of source tree leaves to the right of the new spike.
For step 4, we use
$y$ to mark valid positions for the insertion of subsequent vertices, as we did previously.
The generating function is 
$$
D_\textsf{2}(v) \;=\; \frac{z\+y^2\+u^2}{1-z\+y\+u}\+\frac{D_\textsf{1}(v)-D_\textsf{1}(u)}{v-u} .
$$
Step 3: The effect of adding additional vertices to the right of the spike but to the left of some
other leaves in the bottom subtree is represented by the recursive functional equation
$$
D_\textsf{3}(y,v) \;=\; D_\textsf{2}(v) \:+\: z\+y\+v\+\frac{D_\textsf{3}(y,1)-D_\textsf{3}(y,v)}{1-v} .
$$
Again, the kernel method can be used to solve this for $D_\textsf{3}(y,1)$, the kernel being cancelled by setting $v=1/(1-z\+y)$.

Step 4: Finally, the addition of vertices to the right of the bottom subtree is reflected by the {\L}uka\-sie\-wicz operator $\oper_{\textsf{\L}}$, giving rise to the recursive functional equation
$$
D_\textsf{4}(y) \;=\; D_\textsf{3}(y,1) \:+\: z\+y^2\+\frac{D_\textsf{4}(1)-D_\textsf{4}(y)}{1-y} ,
$$
which can be solved for $D_\textsf{4}(1)$ by cancelling the kernel with $y=(1-\sqrt{1-4\+z})/2\+z$.

The generating function for the set of permutations resulting from the second way of creating a \prefix{$\CCC$}permutation from an \prefix{$\AAA$}permutation
is then defined by
$$
C_{\textsf{AC2}}(u) \;=\; D_\textsf{4}(1) - B_{\textsf{AB2}}(u) .
$$

Our work is almost complete. We only have to consider how a source graph may be added to a \prefix{$\CCC$}permutation.
In fact, the situation is extremely constrained.
First, as noted earlier, the source graph must be positioned to the left of the spike.
Furthermore, the presence of a $\mathbf{13524}$ or $\mathbf{14523}$ means that the addition of a source graph with even a single non-root vertex would create a $\mathbf{1234}$.
So the only possibility is the addition of a trivial (single vertex) source tree.
Thus the bivariate generating function $C(u)$ of set $\CCC$ is defined by the following recursive functional equation:
$$
C(u)  \;=\;  C_\ssS(u) + C_{\textsf{AC1}}(u) + C_{\textsf{AC2}}(u) \:+\: z\+u\+\frac{C(1)-C(u)}{1-u}.
$$
where the final term reflects the addition of a trivial source tree to a \prefix{$\CCC$}permutation.
This equation can be solved to yield the following expression for $C(1)$ by 
a sixth and final application of the kernel method,
cancelling the kernel by setting $u=1/(1-z)$:
$$
\frac{-1+10\+z-35\+z^2+52\+z^3-35\+z^4+12\+z^5 \:+\: (1-8\+z+21\+z^2-22\+z^3+11\+z^4-2\+z^5)\sqrt{1-4\+z}}{2\+z^3\+(1-4\+z)\+(1-3\+z+z^2)} .
$$

We now have all we need to prove Theorem~\ref{thmF} by
obtaining an explicit expression for the generating function
that enumerates class
$\FFF$.
Since $\FFF$ is the disjoint union of $\AAA$, $\BBB$ and $\CCC$,
its generating function
is
given by
$A(1)+B(1)+C(1)$.
Thus, by appropriate expansion and simplification,
the generating function for
$\av(\mathbf{1234},\mathbf{2341})$
can be shown to be equal to
$$
\frac
{2-10\+z+9\+z^2+7\+z^3-4\+z^4 \:-\: (2-8\+z+9\+z^2-3\+z^3)\+\sqrt{1-4\+z}}
{(1-3\+z+z^2)\+\big((1-5\+z+4\+z^2) \:+\: (1-3\+z)\+\sqrt{1-4\+z}\big)}.
$$
This has singularities at $z=\frac{1}{4}$, $z=\frac{1}{2}\+(3-\sqrt{5})$ and $z=\frac{1}{2}\+(3+\sqrt{5})$. Hence, the growth rate of $\av(\mathbf{1234},\mathbf{2341})$ is equal to 4, the reciprocal of the least of these.

The first twelve terms of the sequence $|\FFF_n|$ are 1, 2, 6, 22, 89, 376, 1611, 6901, 29375, 123996, 518971, 2155145.
More values
can be found at
\href{http://oeis.org/A165540}{A165540} in OEIS~\cite{OEIS}.

\vspace{9pt}
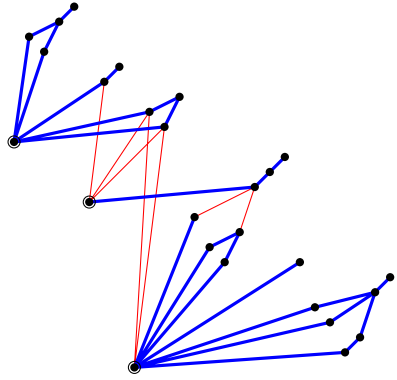
\begin{figure}[ht] 
  $$
  \begin{tikzpicture}[scale=0.20,line join=round]
    \draw [red,thin] (6,15)--(7,23);
    \draw [red,thin] (6,15)--(11,20)--(9,4)--(10,21)--(6,15);
    \draw [red,thin] (13,14)--(17,16)--(16,13);
    \draw [blue,very thick] (8,24)--(7,23)--(1,19)--(2,26)--(4,27)--(3,25)--(1,19);
    \draw [blue,very thick] (4,27)--(5,28);
    \draw [blue,very thick] (1,19)--(10,21)--(12,22)--(11,20)--(1,19);
    \draw [blue,very thick] (6,15)--(17,16)--(19,18);
    \draw [blue,very thick] (13,14)--(9,4)--(14,12)--(16,13)--(15,11)--(9,4);
    \draw [blue,very thick] (26,10)--(25,9)--(21,8)--(9,4)--(22,7)--(25,9)--(24,6)--(23,5)--(9,4)--(20,11);
    \plotpermnobox{}{19,26,25,27,28,15,23,24,4,21,20,22,14,12,11,13,16,17,18,11,8,7,5,6,9,10}
    \draw [thin] (1,19)  circle [radius=0.4];
    \draw [thin] (6,15) circle [radius=0.4];
    \draw [thin] (9,4) circle [radius=0.4];
  \end{tikzpicture}
  $$
  \caption{
           A permutation in class $\EEE$, spanned by three source graphs}
  \label{figEPerm}
\end{figure}
\section{Permutations avoiding \texorpdfstring{$\mathbf{1243}$ and $\mathbf{2314}$}{1243 and 2314}}\label{sectE}

Let us 
use $\EEE$ to denote 
$\av(\mathbf{1243},\mathbf{2314})$.
What can we say about the structure of source graphs in~$\EEE$?
Firstly, since $H_\mathbf{1243} =
\raisebox{-2.5pt}{\begin{tikzpicture}[scale=0.12,line join=round]
  \draw[] (1,1)--(2,2);
  \draw[] (3,4)--(2,2)--(4,3);
  \plotpermnobox{}{1,2,4,3}
\end{tikzpicture}}
$
may not occur as a subgraph, only the root of a
source graph may fork towards the upper right.
Secondly, each source graph in $\EEE$ is \emph{plane}.
This is the case because every non-plane graph
contains a
$H_\mathbf{2143} =
\raisebox{-2.5pt}{\begin{tikzpicture}[scale=0.12,line join=round]
  \draw[] (1,2)--(3,4)--(2,1)--(4,3)--(1,2);
  \plotpermnobox{}{2,1,4,3}
\end{tikzpicture}}
$, and, furthermore, any $\mathbf{2143}$ in a source graph occurs as part of a $\mathbf{13254}$ (where the $\mathbf{1}$ is the root of the source graph). But this is impossible in $\EEE$, since $\mathbf{13254}$
does not avoid $\mathbf{1243}$.

\begin{figure}[ht]
  $$
  \begin{tikzpicture}[scale=0.2,line join=round]
    \draw [blue,very thick] (2,18)--(1,1)--(3,17);
    \draw [blue,very thick] (4,15)--(1,1)--(5,14);
    \draw [blue,very thick] (10,11)--(1,1)--(13,10);
    \draw [blue,very thick] (14,8)--(1,1)--(15,5);
    \draw [blue,very thick] (16,4)--(1,1)--(17,2);
    \draw [ultra thick,black!40!green] (5,14)--(6,16)--(4,15);
    \draw [ultra thick,black!40!green] (6,16)--(7,19)--(3,17);
    \draw [ultra thick,black!40!green] (2,18)--(7,19)--(9,21);
    \draw [ultra thick,black!40!green] (10,11)--(12,13);
    \draw [ultra thick,black!40!green] (17,2)--(18,3)--(19,6)--(16,4);
    \draw [ultra thick,black!40!green] (15,5)--(19,6)--(20,7)--(21,9)--(14,8);
    \plotpermnobox{21}{1,18,17,15,14,16,19,20,21,11,12,13,10,8,5,4,2,3,6,7,9}
    \draw [thin] (1,1)  circle [radius=0.4];
    \draw [thick,black!40!green] (9,21)  circle [radius=0.4];
    \draw [thick,black!40!green] (12,13) circle [radius=0.4];
    \draw [thick,black!40!green] (13,10) circle [radius=0.4];
    \draw [thick,black!40!green] (21,9)  circle [radius=0.4];
  \end{tikzpicture}
  $$
  \caption{A source graph for class $\EEE$, constructed from four u-trees}
  \label{figESource}
\end{figure}
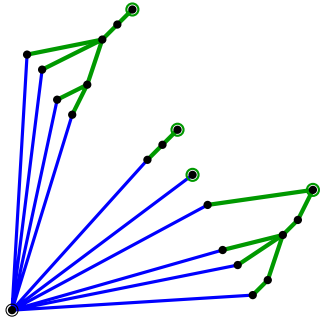
If we combine these two observations, we see that
the non-root vertices of a source graph consist of a sequence of
inverted subtrees whose roots are right-to-left maxima.
The avoidance of $H_\mathbf{2314} =
\raisebox{-2.5pt}{\begin{tikzpicture}[scale=0.12,line join=round]
  \draw[] (1,2)--(2,3)--(4,4)--(3,1);
  \plotpermnobox{}{2,3,1,4}
\end{tikzpicture}}
$
places restrictions on the structure of the subtrees, so that they must consist of a path at the lower right, which we call the \emph{trunk}, with pendant edges attached to its
left.
It is readily seen that these correspond to permutations in
$\av(\mathbf{132},\mathbf{231})$.
We call trees of this form \emph{u-trees}, short for
\emph{unbalanced} trees.
See Figure~\ref{figESource} for an illustration of a source graph constructed from u-trees.

The class $\UUU$ of u-trees satisfies the structural equation
$$
  \UUU \;=\; \ZZZ\times\seq{\ZZZ\times\seq{\ZZZ}}
$$
where the first term on the right represents the lowest leaf at the tip of the trunk and the second represents the remainder of the vertices in the trunk, each with a (possibly empty) sequence of pendant edges attached to the upper left.
Hence the generating function for $\UUU$ is
\begin{equation*}
  U(z) \;=\; \frac{z\+(1-z)}{1-2\+z}.
\end{equation*}
If we use $u$ to mark the number of u-trees, the class $\SSS$ of source graphs
satisfies the structural equation
$$
  \SSS \;=\; \ZZZ\times\seq{u\+\UUU}
$$
and thus has bivariate generating function
\begin{equation*}
  S(u) \;=\; S(z,u) \;=\; \frac{z\+(1-2\+z)}{1-(2+u)\+ z+u\+ z^2}.
\end{equation*}

Let us now examine how a permutation in $\EEE$ can be built from a sequence of source graphs.
Observe that, when a source graph is added,
no vertex of the source graph can be positioned between two vertices of a u-tree in the bottom subgraph, because otherwise
a $\mathbf{2314}$ would be created.
In addition,
there are strong constraints on when u-trees in the new source graph can be positioned to the left of a u-tree in the bottom subgraph.

\begin{figure}[ht]
  $$
  \raisebox{22.5pt}{
  \begin{tikzpicture}[scale=0.2,line join=round]
    \draw [very thick] (2,19)--(1,8)--(3,18);
    \draw [very thick] (5,16)--(1,8)--(6,15);
    \draw [very thick] (9,13)--(1,8)--(10,12);
    \draw [very thick] (12,10)--(1,8)--(13,9);
    \draw [blue,very thick] (15,6)--(8,1)--(16,5);
    \draw [blue,very thick] (18,3)--(8,1)--(19,2);
    \draw [red,thin] (9,13)--(8,1)--(10,12);
    \draw [red,thin] (12,10)--(8,1)--(13,9);
      \plotpermnobox{20}{8,19,18,20,16,15,17,1,13,12,14,10,9,11,6,5,7,3,2,4}
    \draw [very thick,black!40!green,fill=black!40!green!65!white] (3,18)--(4,20)--(2,19);
    \draw [very thick,black!40!green,fill=black!40!green!65!white] (6,15)--(7,17)--(5,16);
    \draw [very thick,black!40!green,fill=black!40!green!65!white] (10,12)--(11,14)--(9,13);
    \draw [very thick,black!40!green,fill=black!40!green!65!white] (13,9)--(14,11)--(12,10);
    \draw [very thick,black!40!green,fill=black!40!green!65!white] (16,5)--(17,7)--(15,6);
    \draw [very thick,black!40!green,fill=black!40!green!65!white] (19,2)--(20,4)--(18,3);
    \draw [thin] (1,8)  circle [radius=0.4];
    \draw [thin] (8,1) circle [radius=0.4];
  \end{tikzpicture}
  }
  \qquad\qquad\qquad\quad
  \begin{tikzpicture}[scale=0.2,line join=round]
    \draw [very thick] (2,24)--(1,13)--(3,23);
    \draw [very thick] (6,21)--(1,13)--(7,20);
    \draw [very thick] (9,18)--(1,13)--(10,17);
    \draw [very thick] (1,13)--(17,14);
    \draw [ultra thick,black!40!green] (19,16)--(17,14);
    \draw [blue,very thick] (12,12)--(5,1)--(13,11);
    \draw [blue,very thick] (14,9)--(5,1)--(15,8);
    \draw [blue,very thick] (20,6)--(5,1)--(21,5);
    \draw [blue,very thick] (23,3)--(5,1)--(24,2);
    \draw [red,thin] (6,21)--(5,1)--(7,20);
    \draw [red,thin] (9,18)--(5,1)--(10,17);
    \draw [red,thin] (12,12)--(17,14)--(13,11);
    \draw [red,thin] (16,10)--(17,14);
      \plotpermnobox{25}{13,24,23,25,1,21,20,22,18,17,19,12,11,9,8,10,14,15,16,6,5,7,3,2,4}
    \draw [very thick,black!40!green,fill=black!40!green!65!white] (3,23)--(4,25)--(2,24);
    \draw [very thick,black!40!green,fill=black!40!green!65!white] (7,20)--(8,22)--(6,21);
    \draw [very thick,black!40!green,fill=black!40!green!65!white] (10,17)--(11,19)--(9,18);
    \draw [very thick,black!40!green,fill=black!40!green!65!white] (15,8)--(16,10)--(14,9);
    \draw [very thick,black!40!green,fill=black!40!green!65!white] (21,5)--(22,7)--(20,6);
    \draw [very thick,black!40!green,fill=black!40!green!65!white] (24,2)--(25,4)--(23,3);
    \draw [thin] (1,13)  circle [radius=0.4];
    \draw [thin] (5,1) circle [radius=0.4];
  \end{tikzpicture}
  $$
  \caption{The two methods for adding a source graph in class $\EEE$; u-trees are shown schematically as filled triangles}
  \label{figEBuild}
\end{figure}
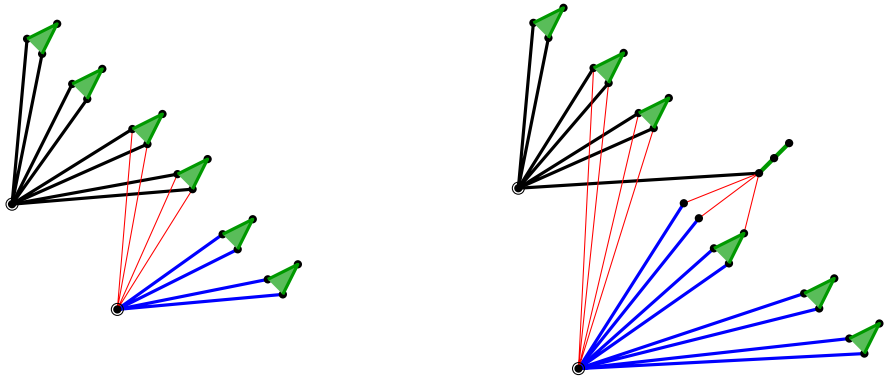
These conditions result in there being two distinct ways in which a source graph may be added.
These are illustrated in Figure~\ref{figEBuild}.
In the first method, the root of the source graph is positioned to the left of zero or more u-trees in the bottom subgraph and the u-trees in the source graph are positioned to the right of the bottom subgraph.

The second method is more subtle. It is only applicable if
the rightmost u-tree of the bottom subgraph
is a path.
If that is the case, then an initial sequence of
u-trees in the source graph can be positioned to the left of this path subtree, as long as
each of them, except possibly the last,
consists of a single vertex.
If
the rightmost u-tree of the bottom subgraph
were not a path, then a $\mathbf{1243}$ would be created.
Similarly, if a non-final u-tree consisted of more than one vertex, then a $\mathbf{2314}$ would be created.

In order to handle this second method, we need to keep track of those source graphs in which the rightmost u-tree is a path.
Let $\SSS_\PP$ be the class of such graphs. It satisfies the structural equation
$$
  \SSS_\PP \;=\; \ZZZ\times\seq{u\+\UUU}\times u\+\seqplus{\ZZZ} ,
$$
where $u$ marks the number of u-trees as before.
This class thus has bivariate generating function
\begin{equation*}
  S_\PP(u) \;=\; S_\PP(z,u) \;=\; \frac{u\+ z^2\+ (1-2\+ z)}{(1-z)\+ \big(1-(2+u)\+ z+u\+ z^2\big)}.
\end{equation*}

In order to distinguish between those situations when the second method of adding a source graph is applicable and those when it isn't,
let us use $\PPP$ to denote the set of those permutations in $\EEE$ whose Hasse graphs have bottom subgraphs in which the rightmost u-tree is a path.

We are interested in determining the two bivariate generating functions
$E(u)=E(z,u)$ and
$P(u)=P(z,u)$ for $\EEE$ and $\PPP$ respectively, where $u$ marks the number of u-trees \emph{in the bottom subgraph}.
To do this, we will establish four linear operators on these generating functions that reflect the different ways in which a source graph can be added.

The action of adding a source graph using the first method is readily seen to be reflected by the following linear operator:
\begin{equation*}
  \oper_{\EE\EE}\big[f(u)\big] \;=\; S(u)\+\frac{f(1)-u\+f(u)}{1-u}.
\end{equation*}
The first method creates a member of $\PPP$ from an arbitrary element of $\EEE$ whenever
the source graph is in $\SSS_\PP$ (i.e.
its rightmost u-tree is a path).
Thus the appropriate linear operator 
is
\begin{equation*}
  \oper_{\EE\PP}\big[f(u)\big] \;=\; S_\PP(u)\+\frac{f(1)-u\+f(u)}{1-u}.
\end{equation*}

Now let us determine the linear operators corresponding to the second method of adding a source graph.

The set, $\SSS^\star$, of source graphs that can be added using the second method
satisfies the structural equation
$$
  \SSS^\star \;=\; \ZZZ\times\seq{\ZZZ}\times u\+\UUU\times\seq{u\+\UUU},
$$
in which the third term on the right identifies the u-tree which is positioned immediately to the left of the rightmost (path) u-tree in the bottom subgraph.
This specification thus counts multiple times those source graphs that can be added in more than one way due to the presence of a non-empty initial sequence of single-vertex u-trees.
Note also that we don't mark the initial sequence of single-vertex u-trees with $u$.
The generating function for $\SSS^\star$ is
\begin{equation*}
  S^\star(u) \;=\; \frac{u\+ z^2}{1-(2+u)\+ z+u\+ z^2}.
\end{equation*}
The action of adding a source graph using the second method is then seen to be reflected
by the following linear operator:
\begin{equation*}
  \oper_{\PP\EE}\big[f_\PP(u)\big] \;=\; S^\star(u)\+\frac{f_\PP(1)-f_\PP(u)}{1-u} .
\end{equation*}
Finally, let us consider when adding a source graph to an arbitrary member of $\PPP$ creates another permutation in $\PPP$.
The second method creates an element of $\PPP$ if the source graph is in $\SSS_\PP$ and its rightmost (path) u-tree is added to the right of the bottom subgraph.
An element of $\PPP$ is also created if the source graph has a single path u-tree or consists of a single vertex (the root).
Thus the set, $\SSS_\PP^\star$, of source graphs, counted with multiplicity,
that can be added to create an element of $\PPP$
satisfies the structural equation
$$
  \SSS_\PP^\star \;=\; \ZZZ\times\seq{\ZZZ}\times\seqplus{u\+\UUU}\times u\+\seqplus{\ZZZ} \:+\: \ZZZ\times u\+\seq{\ZZZ}.
$$
Its generating function is
\begin{equation*}
  S_\PP^\star(u) \;=\; \frac{u\+ z\+ (1-2\+ z)\+(1-u\+z)}{(1-z)\+ \big(1-(2+u)\+ z+u\+ z^2\big)},
\end{equation*}
and the corresponding linear operator
is
\begin{equation*}
  \oper_{\PP\PP}\big[f_\PP(u)\big] \;=\; S_\PP^\star(u)\+\frac{f_\PP(1)-f_\PP(u)}{1-u}.
\end{equation*}

We are now in a position to derive the generating function 
for $\EEE$ and hence prove Theorem~\ref{thmE}.
From the analysis above, we know that
the bivariate generating function $E(u)=E(z,u)$ of class $\EEE$ is defined by the following pair of mutually recursive functional equations:
\begin{equation*}
  \begin{array}{rclcrcr}
   E(u) & = & S(u)     & \!+\! & \oper_{\EE\EE}\big[E(u)\big] & \!+\! & \oper_{\PP\EE}\big[P(u)\big] \\[3pt]
   P(u) & = & S_\PP(u) & \!+\! & \oper_{\EE\PP}\big[E(u)\big] & \!+\! & \oper_{\PP\PP}\big[P(u)\big]
  \end{array}  .
\end{equation*}
These can be expanded to give the following:
\begin{equation}\label{eqnE1}
  E(u) \;=\; z\+\frac{(1-2\+z)\+\big(1-u+E(1)-u\+E(u)\big) \:+\: u\+z\+\big(P(1)-P(u)\big)}{(1-u)\+ \big(1-(2+u)\+ z+u\+ z^2\big)} ,
\end{equation}
\begin{equation}\label{eqnE2}
  P(u) \;=\; u\+z\+(1-2\+z)\frac{z\+\big(1-u+E(1)-u\+E(u)\big)\:+\: (1-u\+z)\+\big(P(1)-P(u)\big)}{(1-u)\+ (1-z)\+ \big(1-(2+u)\+ z+u\+ z^2\big)} .
\end{equation}
An unusual simultaneous double application of the kernel method can then be used
to yield the
algebraic generating function for class $\EEE$ as follows.

First, we eliminate $P(u)$ from \eqref{eqnE1} and \eqref{eqnE2}, and express $E(u)$
in terms of $E(1)$ and $P(1)$ as a rational function.
Cancelling the resulting kernel,
\begin{equation*}
(1-3\+z+2\+z^2)
\:-\:
(2-7\+z+7\+z^2-z^3)\+u
\:+\:
(1-3\+z+3\+z^2)\+u^2
\:-\:
(z-3\+z^2+3\+z^3)\+u^3 ,
\end{equation*}
with the appropriate root
then gives us an equation relating $E(1)$ and $P(1)$.

Secondly, we eliminate $E(u)$ from \eqref{eqnE1} and \eqref{eqnE2}, and express $P(u)$ in terms of $E(1)$ and $P(1)$.
Cancelling the (same) kernel (using a different root)
gives a second equation relating $E(1)$ and $P(1)$.

Finally, we eliminate $P(1)$ from these two equations to yield an extremely complicated explicit expression for $E(1)$. 

Thus, using a computer algebra system to handle the details of the algebraic manipulation, it can be determined that the generating function $F(z)=E(1)$ for 
$\av(\mathbf{1243},\mathbf{2314})$ has the minimal polynomial
$$
      (z-3\+z^2+2\+z^3)
\:-\: (1-5\+z+8\+z^2-5\+z^3)\+F(z)
\:+\: (2\+z-5\+z^2+4\+z^3)\+F(z)^2
\:+\: z^3\+F(z)^3
.
$$
The growth rate of the class is given by the reciprocal of the least positive real singularity of its generating function~\cite[Theorems~IV.6 and~IV.7]{FS2009}. Hence, by determining the location of the singularities of $E(1)$, it is possible to establish that the growth rate of class $\EEE$ is approximately 5.1955, the greatest real root of the quintic polynomial 
$$
2 - 41\+z + 101\+z^2 - 97\+z^3 + 36\+z^4 -4\+z^5 ,
$$
as required.

The first twelve terms of the sequence $|\EEE_n|$ are 1, 2, 6, 22, 88, 367, 1571, 6861, 30468, 137229, 625573, 2881230.
More values
can be found at
\href{http://oeis.org/A165539}{A165539} in OEIS~\cite{OEIS}.

\subsubsection*{Acknowledgements}
Michael Albert's PermLab software~\cite{PermLab} was of particular benefit in helping to visualise and explore the structure of permutations in the two permutation classes.
The author is also grateful to Robert Brignall and two anonymous referees for suggestions that led to improvements in the presentation of parts of the paper.

\emph{S.D.G.}


\bibliographystyle{plain}
{\footnotesize\bibliography{mybib}}

\end{document}